\documentclass[11pt]{amsart}
\usepackage{graphicx}
\usepackage{epsfig}
\usepackage{amsfonts}
\usepackage{amssymb}
\usepackage{amsmath}
\usepackage{amsthm}
\usepackage{latexsym}
\usepackage{amscd}
\usepackage{flafter}
\usepackage{epsf}
\newcommand{\PN}{{\mathbb{P}^N}}
\newcommand{\Group}{\mathfrak{G}}
\newcommand{\GroupH}{\mathfrak{H}}

\newcommand{\ParMod}{\mathcal{P}}
\newcommand{\MapMod}{\mathcal{X}}
\newcommand{\Almost}{\mathcal{J}}
\newcommand{\Inv}{\mathcal{I}}
\newcommand{\lra}{\longrightarrow}
\newcommand{\ra}{\rightarrow}
\newcommand{\Sig}{\Sigma}
\newcommand{\jsig}{{j_\Sigma}}
\newcommand{\Mod}{\mathcal{M}}
\newcommand{\Modbar}{{\overline{\mathcal{M}}}}
\newcommand{\Ci}{{\mathcal{C}}}
\newcommand{\Cibar}{{\overline{\mathcal{C}}}}
\newcommand{\xra}{\xrightarrow}
\newcommand{\Space}{\mathcal{Y}}
\newcommand{\Bundle}{\mathcal{E}}
\newcommand{\Homi}{\text{Hom}^{0,1}}
\newcommand{\dbar}{\overline{\partial}}
\newcommand{\ModX}{\mathcal{M}_g(X,\beta)}
\newcommand{\ModGX}{\mathcal{M}^\Group_g(X,\beta)}

\newcommand{\ModXpar}{\mathcal{M}_g(X,\beta;(J;v))}
\newcommand{\Orbit}{\mathcal{O}}
\newcommand{\Ker}{\text{Ker}}
\newcommand{\Cok}{\text{Coker}}
\newcommand{\m}{\mathfrak{m}}

\newcommand{\DivGm}{\mathcal{D}^\Group_\m}

\newcommand{\DivG}{\mathcal{D}^\Group}
\newcommand{\e}{\epsilon}
\newcommand{\thet}{{\overline{\theta}}}
\newcommand{\muu}{{\overline{\mu}}}

\newcommand{\PathMod}{\mathcal{Q}(\gamma_0,\gamma_1)}
\newcommand{\NodX}{\mathcal{N}_g(X,\beta)}
\newcommand{\NodGX}{\mathcal{N}^\Group_g(X,\beta)}
\newcommand{\NodGXpar}{{\mathcal{N}_g^{\Group}(X,\beta;\gamma)}}
\newcommand{\NodXpar}{\mathcal{N}_g(X,\beta;\gamma)}
\newcommand{\Divl}{\overline{\mathcal{D}}^{\Group}_{\{\m_1,...,\m_\ell\}}}
\newcommand{\Divisor}{\overline{\mathcal{D}}^{\Group}_{\m^0}}

\newcommand\Dual{\mathcal D}
\newcommand\Duality\Dual

\newcommand\relspinc{\underline{\spinc}}

\newcommand\ModSphere{\ModFlow\left({\mathbb S}\longrightarrow
\Sym^{g-1}(\Sigma_{1})\times \Sym^2(\Sigma_{2})\right)}
\newcommand\ModSpheres\ModSphere

\newcommand\UnparModSp{\widehat \ModSp}
\newcommand\UnparModFlow\UnparModSp

\newcommand{\spinc}{\mathfrak s}

\newcommand\ModMaps{\mathcal M}
\newcommand\ModSp\ModMaps

\newcommand\spincrel\relspinc

\newtheorem{thm}{Theorem}[section]

\newtheorem{cor}[thm]{Corollary}

\newtheorem{lem}[thm]{Lemma}

\newtheorem{defn}[thm]{Definition}

\newtheorem{remark}[thm]{Remark}

\def\endproof{\relax\ifmmode\expandafter\endproofmath\else
  \unskip\nobreak\hfil\penalty50\hskip.75em\hbox{}\nobreak\hfil\bull
  {\parfillskip=0pt \finalhyphendemerits=0 \bigbreak}\fi}
\def\endproofmath$${\eqno\bull$$\bigbreak}
\def\bull{\vbox{\hrule\hbox{\vrule\kern3pt\vbox{\kern6pt}\kern3pt\vrule}\hrule}}

\newcommand{\Q}{\mathbb{Q}}
\newcommand{\R}{\mathbb{R}}

\newcommand{\C}{\mathbb{C}}

\newcommand{\Z}{\mathbb{Z}}

\newcommand{\ModSWfour}{\mathcal{M}}
\newcommand{\ModFlow}{\ModSWfour}

\newcommand\Hom{\mathrm{Hom}}

\newcommand\abuts\Rightarrow
\newcommand\Sym{\mathrm{Sym}}

\begin{document}

\title{Embedded curves and Gromov-Witten invariants of three-folds}%
\author{Eaman Eftekhary}%
\address{ Mathematics Department, Harvard University, 1 Oxford Street, Cambridge, MA 02138}%
\email{eaman@math.harvard.edu}%

%\thanks{}%
%\subjclass{}%

\date{September 2004}%
%\dedicatory{}%
%\commby{}%
% ----------------------------------------------------------------
\begin{abstract}
Associated with a prime homology class $\beta \in P_2(X,\Z)$ (i.e.
$\beta=p\alpha$ and $\alpha \in H_2(X,\Z)$ imply $p=1$ or $p$ is
an odd prime) on a symplectic three-manifold with vanishing first
Chern class, we count the embedded perturbed pseudo-holomorphic
curves in $X$ of a fixed genus $g$ to obtain certain integer
valued invariants analogous to Gromov-Witten invariants of $X$. 
\end{abstract}
\maketitle
% ----------------------------------------------------------------
\section{Introduction and main theorems}

The aim of the work in this paper, is to construct some integer
valued invariants of the symplectic threefolds with vanishing
first Chern class, along the lines that Gromov-Witten invariants
are defined by Ruan and Tian (\cite{R-T1,R-T2}). The motivation
for this project is the conjecture of Gopakumar and Vafa, which
writes the Gromov-Witten invariants of a Calabi-Yau threefold $X$
in terms of some (not mathematically defined) integer valued
invariants, called the Gopakumar-Vafa invariants.\\

In their paper \cite{Go-Va-1}, Gopakumar and Vafa introduce
certain counts of the so called BPS-states to get integer valued
invariants of Calabi-Yau threefolds. These are integers
$n_h(\alpha)$ associated with a genus $h$ and a homology
class  $\alpha \in H_2(X,\mathbb{Z})$, called Gopakumar-Vafa invariants.\\

The mathematical definition of these invariants is yet to be
understood. There has been an attempt by Hosono, Saito and
Takahashi \cite{H-S-T} to define these numbers through some
intersection cohomology construction. They are not able to show
the invariance of these numbers and also it is not clear that they
satisfy the Gopakumar-Vafa
equation predicted in \cite{Go-Va-1}.\\

The Gopakumar-Vafa equation is a generating function equation,
which relates Gopakumar-Vafa invariants with Gromov-Witten
invariants. More precisely, it reads as:\\

\begin{equation}\label{eq:GV}
\sum_{\substack{g\geq0\\ \beta \in
H_2(X,\mathbb{Z})}}N_g(\beta)q^\beta
\lambda^{2g-2}= \sum_{k>0} \  \frac{1}{k}\sum_{\substack{h\geq0 \\
\alpha\in H_2(X,\mathbb{Z})
}}n_h(\alpha)(2\text{sin}(\frac{k\lambda}{2}))^{2h-2} q^{k\alpha}.
\end{equation}
\\

One may look at this equation as a definition for  Gopakumar-Vafa
invariants (see \cite{R-Jim}) and then it is still to be shown
that these numbers are in fact
integers.\\

One other way to follow the motivation provided by Gopakumar and
Vafa, is trying to write Gromov-Witten invariants in terms of some
other integer valued invariants of $(X,\omega)$, when $X$ is a
symplectic $3$-fold with vanishing first Chern class, or more, a
Calabi-Yau $3$-fold. For instance, the relation between
Gromov-Witten invariants and Donaldson-Thomas invariants suggested
in \cite{Davesh} follows these lines.\\

Fix a genus $g>1$.  As an attempt toward such a construction, we
introduce the invariants

$$\Inv_g:P_2(X,\Z) \lra \Z,$$
which assign integer numbers to any \emph{prime homology class}
$\beta \in P_2(X,\Z)$. A prime homology class is defined to be any
homology classes $\beta \in H_2(X,\Z)$ such that if $\beta=p
\alpha$ with $1<p\in \Z$ and $\alpha \in H_2(X,\Z)$,
then $p\neq 2$ is a prime number.\\

We will embed the coarse universal curve $\Cibar_g$ of the moduli
space $\Modbar_g$ (of genus $g$ curves) in a projective space $\PN$.
Correspondingly we will consider a generic almost complex
structure $J$ on the tangent bundle of the three-fold $X$ and a
generic perturbation term

$$v\in \Gamma(\PN\times X,p_1^*\Omega_\PN^{0,1}\otimes_Jp_2^*T'X).$$
If the choice of $(J,v)$ is generic enough, there are only
finitely many pairs $(f,j_\Sig)$ of a complex structure $j_\Sig\in
\Mod_g$ on a surface $\Sig$ of genus $g$ and a map $f:\Sig \lra
X$, such that $f_*[\Sig]$ is the fixed prime homology class
$\beta$ and such that

$$\dbar_{j_\Sig,J}f=(\pi_\Sig \times f)^*v.$$
Here $\pi_\Sig :(\Sig,j_{\Sig})\ra \Cibar_g \subset \PN$ is the map from
$(\Sig,j_\Sig)$ to the universal curve. \\

In fact if we denote the space of all possible pairs $(J,v)$, with $J$ an 
almost complex structure  on the tangent bundle of $X$, compatible with $\omega$,
and with $v$ a perturbation term on $\PN\times X$ with respect to $J$, by 
$\ParMod$, then we show

\begin{thm}
\label{thm:definition} Given $\beta \in P_2(X,\Z)$ and $g>1$, there is a
Bair subset $\ParMod_{reg} \subset \ParMod$ such that for any
$(J,v) \in \ParMod_{reg}$, the above problem has only finitely
many solutions. The linearization $L$ of the perturbed
Cauchy-Riemann equation $$\dbar_{\jsig,J}f=(\pi_\Sig\times f)^*v$$
at any solution will have a trivial kernel and a trivial cokernel.
Moreover, to any such solution is assigned a sign
$$\e(f,\jsig)\in \{-1,+1\},$$ coming from the spectral flow $SF$,
which connects $L$ to a complex $\dbar$ operator.\\
\end{thm}
Using the claim of this theorem, we define
\begin{defn}
Suppose  that $\beta \in P_2(X,\Z)$ and that $g>1$ is a fixed
genus. Fix a pair $(J,v)\in \ParMod_{reg}$ as above and define
$$\Inv_g(\beta)=\sum_{(f,\jsig)}\e(f,\jsig),$$ where $\e(f,\jsig)$
is understood to be zero if $(f,\jsig)$ is
not a solution to the above problem
\end{defn}
Furthermore, we show
\begin{thm}
\label{thm:invariance} If $(J_0,v_0)$ and $(J_1,v_1)$ are in
$\ParMod_{reg}$, then the number $\Inv_g(\beta)$ as computed
using $(J_0,v_0)$, is the same as the number $\Inv_g(\beta)$ as
computed using $(J_1,v_1)$. Moreover, this number does not depend
on the the special embedding of $\Cibar_g$ in the specific
projective space $\PN$, and is independent of the
of the choice of the symplectic form $\omega$ in its isotopy class.
\end{thm}

We hope that this construction can be extended to all the homology
classes in $H_2(X,\Z)$. Once this is done, the Gromov-Witten
invariants may be written down in terms of these counts of
embedded curves in  the symplectic three-folds with vanishing
first Chern class, giving
an equation similar to equation~(\ref{eq:GV}) above.\\

{\bf{Acknowledgements.}} I would like to thank Gang Tian for his
continuous  advice and support as my advisor. Many thanks go to
Rahul Pandharipande and Elleny Ionel for very helpful
discussions.\\

%%%%%%%%%%%%%%%%%%%%%%%%%%%%%%%%%%%
%%%%%%%%%%%%%%%%%%%%%%%%%%%%%%%%%%%
%%%%%%%%%%%%%%%%%%%%%%%%%%%%%%%%%%%

\section{A Review of Gromov-Witten Invariants}

Suppose that $\alpha_1,...,\alpha_k,\beta_1,...,\beta_l$ are homology
classes in a symplectic $n$-manifold $X$ and $\alpha\in
H_2(X,\mathbb{Z})$ is a class in the  second homology of $X$ with
$\Z$ coefficients. Let
$$\Phi_{g,\alpha}^X(\alpha_1,...,\alpha_k|\beta_1,...,\beta_l)$$
denote the genus $g$ Gromov-Witten invariant associated with these
elements of $H_*(X,\mathbb{Z})$  as defined by Ruan and Tian in
\cite{R-T1}. In fact, to define these invariants, consider some
fixed pseudo-manifold representatives of the homology classes
$\alpha_i$ and $\beta_j$, which we still denote by $\alpha_i$ and
$\beta_j$, for $i=1,...,k$ and $j=1,...,l$. Also fix a complex
structure $j$ on a surface $\Sigma$ of genus $g$, together with
$k$ marked points $x_1,...,x_k$ on $\Sigma$. Associated with the
symplectic structure $\omega$ on $X$ is the space of compatible
almost complex structures $J$ on the tangent bundle of $X$, which
will be denoted by $\Almost$. With $j$ and $J\in \Almost$ fixed,
by a perturbation term we mean an element $$v\in
\Gamma(\Sigma\times X,p_1^*\Omega_j^{0,1}(\Sig)\otimes_J
p_2^*T'X).$$

Ruan and Tian fix $\Sig,j,x_1,...,x_k$ and choose  a generic pair
$(J,v)$ of an almost complex structure and a corresponding
perturbation term, and count the number of solutions to the
following problem:

\begin{equation}\label{eq:GWdefenition}
\begin{split}
&f:\Sig \lra X,\\
&f_*[\Sig]=\alpha \in H_2(X,\Z),\\
&\dbar_{j,J}f=df+J\circ df\circ j=(Id \times f)^*v, \\
&f(x_i)\in \alpha_i,\ \ \ \ i=1,...,k,
\end{split}
\end{equation}
such that the image $f(\Sig)$ intersects all of $\beta_j$'s. With
an appropriate choice of the sign, this will give an invariant of
the isotopy class of $\omega$ and the homology classes
$\alpha,\alpha_i,\beta_j$, and the genus $g$, if $X$ is a semi-positive
symplectic manifold and

\begin{displaymath}
\sum_{i=1}^k (2n-\text{deg}(\alpha_i)) +
\sum_{j=1}^{l}(2n-2-\text{deg}(\beta_i))
=2c_1(X)\alpha +2n(1-g).
\end{displaymath}

Using these invariants, they define the quantum cup product and show
that it is associative (see \cite{R-T1}).\\

Later, they defined the higher genus Gromov-Witten  invariants
coupled with gravity. The main idea is to relax the almost complex
structure $j$ on the punctured Riemann surface
$$\Sig-\{x_1,x_2,...,x_k\}$$
in the equation~(\ref{eq:GWdefenition}) to vary in the moduli
space $\Mod_{g,k}$ of $k$-pointed Riemann surfaces of genus $g$.
This way, the classes $\alpha_i$ and $\beta_j$ will be treated in
a uniform way,  so we may assume that $l=0$. The result is an
invariant

$$\Psi^X_{\alpha,g,k}:(H_*(X,\Z))^k\lra \Q,$$
which is zero unless the dimension criterion is satisfied:

\begin{displaymath}
\sum_{i=1}^k (2n-2-\text{deg}(\alpha_i))=2c_1(X)\alpha +2(n-3)(1-g).
\end{displaymath}

 These are
rational numbers because of the orbifold structure on the
universal curve $\Ci_{g,k}$ of the moduli space
$\mathcal{M}_{g,k}$ of genus $g$
Riemann surfaces with $k$ marked points.\\

 In fact,  in the definition of
Gromov-Witten invariants coupled with gravity by Ruan and Tian \cite{R-T2},
they consider the compactified moduli space of $k$ pointed genus $g$ Riemann
surfaces $\overline{\mathcal{M}}_{g,k}$ and the universal curve

\begin{displaymath}
\pi_{g,k}:\overline{\mathcal{C}}_{g,k}\longrightarrow \overline{\mathcal{M}}_{g,k}
\end{displaymath}
over it. The moduli spaces
$\overline{\mathcal{M}}_{g,k},\overline{\mathcal{C}}_{g,k}$ are
both projective varieties sitting in a projective space which we
denote by $\mathbb{P}^N$. A natural candidate for the perturbation
terms associated with an almost complex structure $J$ on $X$ can
be an element of

$$\Hom_J^{0,1}(T\mathbb{P}^N,TX)=
\Gamma(\PN\times X,p_1^*\Omega^{0,1}_{\PN}\otimes_J p_2^*T'X).$$

Then, one will consider the maps $f$ from a Riemann surface
$(\Sigma,j_\Sigma,x_1,...,x_k)\in \overline{\mathcal{M}}_{g,k}$ to
$X$ satisfying $f_*[\Sigma]=\alpha\in H_2(X,\mathbb{Z})$, with certain
constrains $f(x_i)\in \alpha_i$ as above, which also satisfy

\begin{displaymath}
\bar \partial_{J,j_\Sigma}f=(\pi_\Sigma \times f)^*v.
\end{displaymath}
Here $\pi_\Sigma:(\Sigma,j_\Sigma,x_1,...,x_k) \rightarrow \overline{\mathcal{C}}_{g,k}$ is the map
to the universal curve.\\

However, these perturbation terms do not rule out the technical difficulties in the construction
of Ruan and Tian (\cite{R-T2}). The problem is that if a Riemann surface with
marked points $(\Sigma,j_\Sigma;x_1,...,x_k)$ admits an automorphism group $\mathfrak{G}$,
then the map $\pi_\Sigma$ will factor through $\frac{\Sigma}{\mathfrak{G}}$ and the transversality
and compactness arguments fail to work.\\

One is forced to go to a finite "fine" cover of the moduli space
 $\overline{\mathcal{C}}_{g,k}$ to achieve the transversality and compactness. Naturally, one has to
 divide the final count of pseudo-holomorphic curves by the degree of this covering. This will cook up
 invariants of the symplectic manifold which are rational numbers.\\

In particular for a Calabi-Yau threefold, the only chance for
$\Psi$ to be non-zero is for $k=0$. Here associated with each
$\alpha\in H_2(X,\mathbb{Z})$ and each genus $g$, we will get a number
$N_g(\alpha)=\Psi^X_{\alpha,g,0}()$ which is naturally a rational number
as noted above.\\

In this paper, we will try to avoid this detour to fine covers of
the moduli space $\Cibar_g$. The special nature of dimension three
and the vanishing of the first Chern class $c_1(X)$ will help us
out of some of the technical difficulties, at least for the prime
homology classes.\\

The main idea  we want to pursue in this paper  comes from the
work of Taubes in relating Seiberg-Witten invariants to some
Gromov type invariants, which we will discuss now.\\

Suppose that $X$ is a symplectic manifold of real dimension $4$
with a symplectic form $\omega$. Taubes has defined in
\cite{Taubes}, an invariant
\begin{equation}
\text{Gr}: H_2(X,\Z) \lra \Z
\end{equation}
which assigns an integer to any homology class in $H_2(X,\Z)$.
The definition of these invariants relies on counting the
pseudo-holomorphic submanifolds of $X$ with respect to a generic
almost complex structure $J$, which is compatible with $\omega$.
The special nature of dimension $4$ is used to exclude the
convergence of a sequence of pseudo-holomorphic submanifolds to a
singular pseudo-holomorphic curve, or a convergence of the curves
in a class $n\beta$ to a curve, multiply covered by this sequence,
and representing the class $\beta \in H_2(X,\Z)$. The special case
of  the convergence of tori, with topologically trivial normal
bundle, and representing a class $n\beta$, to a torus in class
$\beta$, is more delicate to exclude. As is clear from the above
discussion, these invariants are pretty much similar to the
Gromov-Witten invariants
of Ruan and Tian (the "coupled with gravity" version).\\

Taubes assigns certain weights to each of these tori, which may be
different from $\pm 1$. These signs enable him to conclude that a
passage from an almost complex structure $J_0$, which is generic
enough so that there are only finitely many pseudo-holomorphic
curves in a class $\beta \in H_2(X,\Z)$ with respect to $J_0$, to
another generic almost complex structure
$J_1$, does not change the total count.\\

These invariants are shown to be equivalent to the Seiberg-Witten
invariants for manifolds with $b_2^+>1$ in the celebrated papers
\cite{Taubes1,Taubes2} of Taubes, and the non triviality of them
is a result of the non triviality of the Seiberg-Witten
invariants.\\

We will  start a similar project in dimension $6$ (complex
dimension $3$) through the next couple of sections.\\

%%%%%%%%%%%%%%%%%%%%%%%%%%%%%%%%%%%
%%%%%%%%%%%%%%%%%%%%%%%%%%%%%%%%%%%
%%%%%%%%%%%%%%%%%%%%%%%%%%%%%%%%%%%
\section{Construction of the invariants}

In this section we will introduce an  invariant$$\Inv_g:P_2(X,\Z)
\lra \Z$$ for any genus $g>1$, which assigns an integer to any
\emph{prime} homology class $\beta\in P_2(X,\Z)$. A prime homology
class is defined to be any homology classes $\beta \in H_2(X,\Z)$
such that if $\beta=p \alpha$ with $1<p\in \Z$ and $\alpha \in
H_2(X,\Z)$, then $p\neq 2$ is a prime number. We hope that this
construction can be extended to all the homology classes
in $H_2(X,\Z)$.\\

Let  us begin by embedding the moduli space of genus $g$ Riemann
surfaces with one marked point in some projective space $\PN$.
More precisely, suppose that $\Mod_g$ and $\Modbar_g$ denote the
moduli space of genus $g$ curves and its compactification. The
coarse universal curve over $\Mod_g$ and $\Modbar_g$ will be
denoted by $\Ci_g$ and $\Cibar_g$, respectively. It is a
well-known fact that these moduli spaces are in fact projective
varieties. Thus, we may assume that $\Cibar_g$ is embedded in some
projective space $\PN$. Fix such an embedding. For a complex curve
$(\Sig,\jsig)$, note that there is a natural holomorphic map
$\pi_\Sig$ defined as

$$\pi_\Sig=q_\Sig \circ \tau_\Sig: (\Sig,\jsig) \xra{\tau_\Sig}
\frac{(\Sig,\jsig)}{\text{Aut}(\Sig,\jsig)}
\xra{q_\Sig} \Cibar_g \subset \PN,$$
where $\tau_\Sig$ is the map going from $\Sig$ to its quotient
by its automorphism group and $q_\Sig$ is the embedding
of the quotient in the universal curve $\Cibar_g$.\\

If $(X,\omega)$ is a symplectic threefold with vanishing first
Chern class, i.e. $c_1(X)=0$, we may denote the space of almost
complex structures, compatible with $\omega$, by
$\Almost=\Almost_\omega  $. The elements of $\Almost_\omega$ are
the homomorphism $$J:TX \lra TX$$ with the property that $J\circ
J=- \text{Id}$ and $\omega(\zeta,J\zeta)>0$
for any $\zeta \in T_xX$.\\

Given any compatible almost complex structure $J\in \Almost$, we
may consider the space of associated perturbation terms

$$\Hom_J^{0,1}(T\PN,TX)=\Gamma(\PN \times X,p_1^*\Omega_{\PN}^{0,1}
\otimes_J p_2^*T'X),$$
where $p_1,p_2$ denote the projections to the first and the
second factor, respectively.\\

For a fixed genus $g>1$ and a prime homology class $\beta \in
P_2(X,\Z)$, we are interested in counting the solutions to the
following enumerative problem: Find the number of all somewhere
injective maps $f$ such that

\begin{equation}\label{eq:Cauchy-Riemann-eq}
\begin{split}
&f:(\Sig,\jsig) \lra X\\
&f_*[\Sig]=\beta \in P_2(X,\Z)\\
&\dbar_{\jsig,J}f=(\pi_\Sig\times f)^*v,
\end{split}
\end{equation}
where $J$ is a fixed generic almost complex structure on $X$,
compatible with $\omega$ and $v$ is a fixed
generic perturbation term $v\in \Hom_J^{0,1}(T\PN,TX)$.\\

Let us denote the space of all possible pairs $(J,v)$ by
$\ParMod$. One may think of $\ParMod$ as an infinite dimensional
bundle over $\Almost$. The first major theorem of this chapter is:

\begin{thm}
\label{thm:definition} Given $\beta \in P_2(X,\Z)$, there is a
Bair subset $\ParMod_{reg} \subset \ParMod$ such that for any
$(J,v) \in \ParMod_{reg}$, the above problem has only finitely
many solutions. The linearization $L$ of the perturbed
Cauchy-Riemann equation $$\dbar_{\jsig,J}f=(\pi_\Sig\times f)^*v$$
at any solution will have a trivial kernel and a trivial cokernel.
Moreover, to any such solution is assigned a sign
$$\e(f,\jsig)\in \{-1,+1\},$$ coming from the spectral flow $SF$,
which connects $L$ to a complex $\dbar$ operator.\\
\end{thm}

Using the claim of this theorem, we define
\begin{defn}
Suppose  that $\beta \in P_2(X,\Z)$ and that $g>1$ is a fixed
genus. Fix a pair $(J,v)\in \ParMod_{reg}$ as above and define
$$\Inv_g(X,\beta)=\sum_{(f,\jsig)}\e(f,\jsig),$$ where $\e(f,\jsig)$
is understood to be zero if $(f,\jsig)$ is
not a solution to the above problem
\end{defn}

The second major result is the following:

\begin{thm}
\label{thm:invariance} If $(J_0,v_0)$ and $(J_1,v_1)$ are in
$\ParMod_{reg}$, then the number $\Inv_g(X,\beta)$ as computed
using $(J_0,v_0)$, is the same as the number $\Inv_g(X,\beta)$ as
computed using $(J_1,v_1)$. Moreover, this number does not depend
on the the special embedding of $\Cibar_g$ in the specific
projective space $\PN$, and is independent of the choice of $\omega$ in
its  isotopy class.
\end{thm}

Once we can prove the first claim, the other claims are completely
standard in the theory of Gromov-Witten invariants. Thus we will
not come back to this issue.

%%%%%%%%%%%%%%%%%%%%%%%%%%%%%%%%%%%
%%%%%%%%%%%%%%%%%%%%%%%%%%%%%%%%%%%
%%%%%%%%%%%%%%%%%%%%%%%%%%%%%%%%%%%

\section{Transversality}

In this section we will formulate the above moduli problem in
terms of zeros of a certain section of an infinite dimensional
bundle over some Banach space. In fact, the spaces we will be
dealing with are not Banach spaces, but one can take Sobolov
completions of  these spaces and prove regularity results, as is
standard in the theory of Gromov-Witten invariants (c.f.
\cite{R-T1,R-T2}). We choose to hide this detour to Banach spaces,
in order to simplify the descriptions.\\

To begin, note that the moduli spaces $\Cibar_g$ and $\Modbar_g$
are stratified according to the automorphism group of the curves
$(\Sig,\jsig)\in \Modbar_g$. For a fixed topological action of a
group on a (possibly nodal) Riemann surface $\Sig$, denoted by
$\Group$, let $\Mod_g^\Group$ denote the (open) subvariety of
$\Modbar_g$ , consisting of the curves $(\Sig,\jsig)$ with
automorphism group $\Group$, with the topological type of the
action being specified by $\Group$.\\

Denote the part of the coarse moduli space $\Cibar_g$ that
lies over $\Mod_g^\Group$ by $\Ci_g^\Group$.\\

Let $\MapMod$ denote the space of maps $f$ from $\Sig$ to $X$,
which are somewhere injective and represent the prime homology
class $\beta \in P_2(X,\Z)$. Fix the topological action $\Group$
and denote the topological quotient by $C=\frac{\Sig}{\Group}$.
Denote by $\MapMod^\Group$ the space of maps representing the
class $\beta$ which factor through the quotient map $\tau_\Sig$:

$$f=g \circ \tau_\Sig:\Sig \xra{\tau_\Sig} C=\frac{\Sig}{\Group}\xra{g} X, $$
with $g$ somewhere injective. In particular, the domain of any map
in $\MapMod^\Group$ will be contained in the closure of
$\Ci_g^\Group$. It is then implied that $\beta=|\Group|.\alpha$
for some class $\alpha \in H_2(X,\Z)$. So, either the action is
trivial, or the underlying group of the action is
$\Z_p=\frac{\Z}{p\Z}$ for a prime number $p\neq 2$. We will keep
denoting this action by $\Group$, except if it is necessary to do
otherwise, for two reasons. The first one is to distinguish the
underlying group, from the topological action, which contains more
combinatorial information. The second reason is that, most of the
work in this paper does not use the fact that this group is $\Z_p$,
and may be used in a more general context.\\

We make a remark that for the moduli space $\MapMod^\Group$, we
will also include the maps of the form

$$f=g\circ \tau:(\Sig,\jsig)\xra{\tau} \frac{\Sig}{\Group}\xra{g} X$$
where $\Group$ is a quotient of the automorphism group $\GroupH$
of $(\Sig,\jsig)$, but the map $g$  can not be factored through
any other nontrivial quotient of $\frac{\Sig}{\Group}$. Moreover
the topological action of $\Group$ on $(\Sig,\jsig)$ is assumed
implicit in the notation.\\

Over the appropriate Banach versions of these moduli spaces of
maps (consisting of $(k,p)$-maps) one may construct a Banach
bundle as follows: First of all let $\Almost$ denote the space of
almost complex structures on $X$. Let $\ParMod$ be the bundle over
$\Almost$ defined over a point $J\in \Almost$ by

$$\ParMod_J=\Homi_J(T\PN,TX).$$
Let $\Space^\Group$ be the space $\ParMod\times \Mod_g^\Group
\times \MapMod^\Group$. For a point

$$\mu=((J;v),\jsig, (f:\Sig
\xra{\tau_\Sig} \frac{\Sig}{\Group}\xra{g} X))\in \Space^\Group$$
define the bundle $\Bundle^\Group$ by setting

$$[\Bundle^\Group]_\mu:=
\Gamma_\Group(\Sig,\Omega^{0,1}_{\jsig} \otimes f^*T'X),$$ where
$\Gamma_\Group$ denotes the $\Group$-invariant sections of the
corresponding bundle. In reality, $\Bundle$ should be chosen to be
the bundle of sections in an appropriate Sobolov space. More precisely,
we should consider $(k,p)$-sections of these bundles and then prove
some regularity results. The details of these arguments will be left
to the reader.\\

The Cauchy-Riemann operator defines a section of this bundle by setting:

\begin{equation}\label{eq:Cauchy-Riemann }
\begin{split}
&\dbar_\Group: \Space^\Group \lra \Bundle^\Group,\\
&\dbar_\Group(\mu)=\dbar_{\jsig,J}f-(\pi_\Sig\times f)^*v, \end{split}
\end{equation}
where $\mu=((J;v),\jsig, (f:\Sig \xra{\tau_\Sig}
\frac{\Sig}{\Group}\xra{g} X))$ is the above point of
$\Space^\Group$. It is easy to use the fact that $g$ is somewhere
injective and that the target bundle is the space of
$\Group$-invariant sections to show that the section
$\dbar_\Group$ is everywhere transverse to the zero section. The
proof of this fact is identical to the transversality arguments of
\cite{R-T1,R-T2} and
we choose not to repeat them here.\\

Denote the intersection of $\dbar_\Group$ with the zero section by
$\ModGX$. There is a projection map from $\ModGX$ to the parameter
space $\ParMod$ which we will denote by $\text{proj}_\Group$. This
projection map is a Fredholm operator on the Banach versions of
these moduli spaces. To compute the index, note that the kernel
and the cokernel of the linearization of $\text{proj}_\Group$ are
isomorphic to the kernel and the cokernel of the the linearization
of $\dbar_\Group$ restricted to a fiber of $\text{proj}_\Group$.
This index computation is a special case of the index computation
in the Appendix. It is implied that the index of
$\text{proj}_\Group$ is zero. As a result, for a generic choice of
the parameters $(J,v)\in \ParMod$,
$\text{proj}_\Group^{-1}\{(J,v)\}$
will be a set of isolated solutions to the Cauchy-Riemann
equation~\ref{eq:Cauchy-Riemann-eq }.\\

The parameters $(J,v)$ which are generic in the above sense are
not yet generic enough to rule out the possibility that a sequence
of embedded solutions converges to a multiply covered one. Namely,
a sequence in $\text{proj}_\Group^{-1}\{(J,v)\}$ can converge to a
curve in $\text{proj}_\GroupH^{-1}\{(J,v)\}$, if $\Group$ is a
normal subgroup of $\GroupH$. In particular, although the
solutions in each stratum are isolated, there can still exist
infinitely many of them. A careful study of the behavior of
sequences in $\text{proj}_\Group^{-1}\{(J,v)\}$ is the next step toward
the definition of the embedded  curve invariants.\\

%%%%%%%%%%%%%%%%%%%%%%%%%%%%%%%%%%%
%%%%%%%%%%%%%%%%%%%%%%%%%%%%%%%%%%%
%%%%%%%%%%%%%%%%%%%%%%%%%%%%%%%%%%%
\section{Compactness }

In this section we will fix an almost complex structure $J$ on
$TX$ and a perturbation term $v \in \Homi_J(T\PN,TX)$, without any
genericity assumptions. The problem we want to study, is the
behavior and the possible limits of a sequence of maps
$f_n:(\Sig,j_n)\ra (X,J)$ which satisfy the equations

\begin{displaymath}\begin{split}
&\dbar_{j_n,J}f_n=(\pi_n\times f_n)^*v,\\  \text{where}
\ \ \ &\pi_n:(\Sig,j_n) \lra \Cibar_g \subset \PN.
\end{split}
\end{displaymath}
As usual we will assume that $(f_n)_*[\Sig]$ is a fixed homology
class $\beta \in H_2(X,\Z)$.\\

By Gromov compactness theorem, there is a subsequence of
$\{f_n\}$, which we will identify by the sequence itself, which
converges to a map $$f=f_\infty:(\overline{\Sig},j)\lra X$$
satisfying
$$\dbar_{j,J}f=(\pi\times f)^*v,\ \ \ \ \ \ \
\pi:(\overline{\Sig},j) \lra \Cibar_g \subset \PN.$$

Here $\overline{\Sig}$ can lie in the boundary components of the
moduli space $\Modbar_g$ consisting of nodal curves.\\

If we fix the topological type of the nodal curve
$\overline{\Sig}=\{\Sig^1,...,\Sig^n\}$, where $\Sig^i$ are the
irreducible components of $\overline{\Sig}$, and the automorphism
group $\overline{\Group}$ of the domain, then the moduli of maps
$f=f_\infty:(\overline{\Sig},j)\lra X$ as above, which solve the
perturbed Cauchy-Riemann equation, may be described as the
intersection of a certain $\dbar$ section, and the zero section of
a bundle, similar to the construction of the previous section. The
index of the projection map will be $2-2n$. This implies that for
a generic choice of $(J,v)$, we can not have a solution with a
singular domain. Moreover, if $(J_t,v_t)$ is a generic choice of a
path of parameters (with the fixed generic end points) then we may
assume that the singular solutions do not exist for any of
$(J_t,v_t)$. For a more careful discussion of these boundary
components we refer the reader to \cite{R-T1,R-T2}.\\

The above discussion shows that we only need to worry about the
convergence of a sequence $\{f_n\}$, as above, to a solution
$f=f_\infty:({\Sig},j)\lra X$ with $j_n \ra j$ as $n$ goes to
infinity. If the map  $f=f_\infty:({\Sig},j)\ra X$, is somewhere
injective, the usual transversality results will work to find a
non trivial element in the kernel of the linearization of the
Cauchy-Riemann operator. More precisely, if

$$\mu_n=((J;v),j_n, ({\Sig}\xra{f_n} X))\in \ModX,$$ where
$\ModX=\ModGX$ for $\Group=\{\text{Id}\}$, then $\mu_n$ will
converge to
$$\mu=((J;v),j, ({\Sig}\xra{f} X)).$$ This $\mu$ will also be in
$\ModX$ if $f$ does not factor through a nontrivial quotient map
$\tau:\Sig \ra \frac{\Sig}{\Group}$. This is a contradiction if
$(J,v)$ is generic in the sense of previous section, so that the
elements of proj$_{Id}^{-1}(J,v)$
are isolated. \\

If a sequence $\mu_n=((J;v),j_n,
({\Sig}\xra{f_n} X))\in \ModX$ converges to some $\mu=((J;v),j,
({\Sig}\xra{f} X))$ which is not in $\ModX$, then $\mu$ will be in
$\ModGX$ for some group $\Group$, which is a normal subgroup of
the automorphism group of $(\Sig,j)$. Thus, $f$ will be decomposed
as
$$f:\Sig \xra{\tau_\Sig} \frac{\Sig}{\Group}=C\xra{g}X.$$

After choosing a connection, for large values of $n$, we may write
$(f_n,j_n)$ as the image of some elements

$$(\zeta_n\oplus\eta_n)\in \Gamma(\Sigma,f^*TX) \oplus H^1(\Sigma,T\Sigma)$$
under the exponential map

$$\text{exp}:\Gamma(\Sigma,f^*TX)
\oplus H^1(\Sigma,T\Sigma) \lra \MapMod \times \Mod_g.$$

Note that the linearization of the map

\begin{displaymath}
\begin{split}
&\dbar_{J,v}:\Mod_g \times \MapMod \lra \Bundle,\\
&\dbar_{J,v}(j_\Sig,f)=\dbar_{j_\Sig,J}f-(\pi_\Sig \times f)^*v
\end{split}
\end{displaymath}
may be composed with the projection of

$$T_{(\nu;0)}\Bundle=\Bundle_\nu \oplus T_\nu
(\Mod_g \times \MapMod), \ \ \ \ \nu=(j_\Sig,(f:\Sig \ra X))$$
to the first component to give a map

\begin{displaymath}
\begin{split}
d\dbar_{J,v}:T_\nu(\Mod_g \times \MapMod)
=\Gamma(\Sigma,f^*TX) \oplus &H^1_{\jsig}(\Sigma,T\Sigma)\\
&\lra \Bundle_\nu=\Gamma(\Sig,\Omega^{0,1}_\Sig\otimes f^*T'X).
\end{split}
\end{displaymath}

This decomposition of $\Bundle_\nu$ can be made since
we are looking at a zero of the map $\dbar_{J,v}$.\\

The map $d\dbar_{J,v}$ is a Fredholm operator. As a result, it has
a finite dimensional kernel and a finite dimensional cokernel. In
fact, the index of $d\dbar_{J,v}$ is zero, which says that the
kernel and the cokernel have  the same dimension.\\

Since the sections $(\zeta_n\oplus\eta_n)$ are in the domain of
$d\dbar_{J,v}$, we may look at the projection of these sections to
the kernel of $D=d\dbar_{J,v}$. Let us write

$$(\zeta_n\oplus \eta_n)=(\zeta_n^0\oplus \eta_n^0)+
(\zeta_n^1\oplus \eta_n^1),$$ with $(\zeta_n^0\oplus \eta_n^0)$
the projection on the kernel of $D$, and $(\zeta_n^1\oplus
\eta_n^1)$ in the image of the adjoint operator $D^*$. The
following lemma reproves the non triviality of the kernel of $D$:

\begin{lem}
\label{lem:uniqueness}
 Suppose that $(f_i,j_i)=\text{\emph{exp}}_{(f,j)}(\zeta_i\oplus \eta_i)
 =\text{\emph{exp}}_{(f,j)}(\theta_i)$ for  $i=1,2$,
 be two elements of $\MapMod \times \Mod_g$ such that
$\dbar_{j_i,J} f_i=(\pi_i \times f_i)^*v$, where $\pi_i:(\Sig,j_i)
\ra \Cibar_g \subset \PN$. Furthermore, assume that
$(\zeta_1-\zeta_2)\oplus (\eta_1-\eta_2)$ is in the image of
$D^*$. Then $\theta_1=\theta_2$ if $\|\zeta_i\|$ and $\|\eta_i\|$
are assumed to be small enough.\\
\end{lem}

\begin{proof} Define the function

$$
\mathcal{F}: \Gamma(\Sigma,f^*TX) \oplus H^1(\Sig,T\Sig)
\lra \Gamma(\Sig,\Omega^{0,1}_\Sig\otimes f^*T'X)
$$
for $\theta=(\zeta,\eta)$ by the formula

\begin{equation}
\mathcal{F}(\theta)=\Phi_\theta(\dbar_{j,J}(\text{exp}_{f}(\zeta))-
(\pi_{\text{exp}_j(\eta)}\times \text{exp}_f(\zeta))^*v),
\end{equation}
where $\Phi_\theta$ is the parallel transport map along the curve

$$\text{exp}_{(f,j)}(t\theta)=(\text{exp}_f(t\zeta),\text{exp}_j(t\eta))$$
for $t\in [0,1]$. A simple computation shows that the differential of
$\mathcal{F}$ at zero is $d\mathcal{F}(0)=D$.\\

Assume that $\theta_1=Q\gamma +\theta_2=\hat \theta+\theta_2$. Here
$$
Q:\text{Im}(D) \lra \text{Im}(D^*)
$$
is a right inverse for $D$. It is clearly
possible to write the difference of $\theta_1,\theta_2$ in this form since
$Q$ is surjective.\\
Since $f_i$ are $(J,v,j_i)$-holomorphic, we know that
$\mathcal{F}(\theta_i)=0,\ i=1,2$. Suppose that $\|Q\|=c_0$. Then:

\begin{equation}
\begin{split}
\|\hat \theta\| &=\|Q\gamma\| \\
               & \leq c_0 \|\gamma\| \\
               &=c_0 \|D \hat \theta \| \\
               &=c_0 \|d\mathcal{F}(0)\hat \theta\|\\
               &=c_0 \|\mathcal{F}(\theta_1)-\mathcal{F}(\theta_2)-
                    d\mathcal{F}(0)(\theta_1 -\theta_2) \|\\
               &\leq c_0\|\mathcal{F}(\theta_1)-\mathcal{F}(\theta_2)-
                    d\mathcal{F}(\theta_2)(\theta_1 -\theta_2)\| +
                    c_0\| d\mathcal{F}(\theta_2)-d\mathcal{F}(0)\|
                    \|\hat \theta\| \\
               &\leq c_0c_1\|\hat \theta\|_{L^\infty}.\|\hat \theta\|
                 +c_0c_2\|\theta_2\|.\|\hat \theta\|.
\end{split}
\end{equation}

Here the constants only depend on $(f,j)$ and on $(J,v)$. If
$\|\hat \theta\|_{L^\infty}$ and $\|\theta_2\|$ are small enough
then this implies that $\|\hat \theta\|$ is zero, or equivalently,
$\theta_1=\theta_2$. This completes the proof.
\end{proof}

As noted before, the first consequence of this lemma is that if we
have a sequence of maps

$$\mu_n=((J;v),j_n, ({\Sig}\xra{f_n} X))\in \ModX$$ for some $(J,v)$, converging
to a multiply covered one

$$\mu=((J;v),j,({\Sig}\xra{\tau}\frac{\Sig}{\Group}\xra{g} X))\in \ModGX,$$
then we obtain a nonzero element in the kernel of the linearization
$d\dbar_{J,v}$ (note that this is different from the
linearization of the equivariant section $\dbar_\Group$).\\

This completes our first step toward understanding the convergence of sequences in
$\ModX$.

%%%%%%%%%%%%%%%%%%%%%%%%%%%%%%%%%%%
%%%%%%%%%%%%%%%%%%%%%%%%%%%%%%%%%%%
%%%%%%%%%%%%%%%%%%%%%%%%%%%%%%%%%%%

\section{A structure theorem for parameters}

In this section we will consider some local models that describe
the set of acceptable pairs $(J,v)$ of almost complex structures
$J$ on $X$ and perturbation terms $v$ associated with $J$. To
begin, fix a  real rank-$6$ bundle $E$ on the Riemann surface
$\Sig$ with $c_1(E)=0$, a topological action of a group $\Group$
on $\Sig$ giving a quotient $C=\frac{\Sig}{\Group}$, and an
isomorphism class of an almost complex structure on $E$. The
bundle $E\ra \Sig \xra{\pi} C$ will be assumed to be of the form
$\pi^*F$, where $F\ra C$ is a bundle over $C$ with trivial first
Chern class. Note that the map $\pi:\Sig \ra C$ is obtained as a
result of fixing the action of
$\Group$ on $\Sig$.  \\

Let $\Almost^\Group_E$ denote the space of all almost complex
structures $J:E\ra E$, which are invariant under the action of
$\Group$, or saying in a different way, the space of almost
complex structures on $F$. Consider the following bundle over
$\Mod^\Group_g\times \Almost^\Group_E$- where $\Mod^\Group_g$ is
the part of $\Mod_g$ consisting of the complex structures on
$\Sig$ which have $\Group$ as a  subgroup of their
automorphism group:  Let $\mathcal{G}^\Group$ be the bundle given
by

$$\mathcal{G}_{(j_\Sig,J)}^\Group=\Gamma_\Group(\Sig \times E,p_1^*
\Omega_{\jsig}^{0,1}\otimes_J p_2^*E)$$ for any complex structure
$j_\Sig$ on $\Sig$ and any almost complex structure $J$ on $F$.
Here $E$ is considered as a bundle over itself.
We will identify some "bad locus variety" in the total space
$\mathcal{G}^\Group$.\\

A point in $\mathcal{G}^\Group$ will be of the form $(\jsig,J;v)$
where $\jsig$ is a complex structure on $\Sig$ with automorphism
group $\Group$ such that the map to the quotient $(\Sig,\jsig)\ra
\frac{\Sig}{\text{Aut}(\Sig)}$ is in fact the fixed map $\pi$, $J$
is a complex structure on $F$ (inducing one on $E=\pi^*F$) and $v$
is a $\Group$-invariant perturbation term as above. Associated
with this data we will look at the linearization operator

$$L=L(\jsig,J;v):\Gamma(\Sig,E)\oplus H^1(\Sig,T\Sig)
\lra \Gamma(\Sig,\Omega^{0,1}_{(\Sig,\jsig)}\otimes_J E)$$
defined for a section $\omega$ of $T\Sigma$ to be

\begin{equation}
\begin{split}
L(\zeta,\eta)(\omega)=\nabla_\omega \zeta&+
J\nabla_{\jsig \omega}\zeta+(J\circ d\pi \circ \eta) (\omega)\\
&+\frac{1}{2}\{(\nabla_\zeta J)(d\pi \circ \jsig) -(\nabla_{(\zeta
\oplus \eta)}v)\}.
\end{split}
\end{equation}

This is a Fredholm operator of index zero. We are interested in
excluding the locus of $(\jsig,J;v)$  where
the operator has a nontrivial kernel (c.f. computation of the Appendix).\\

If $\theta=(\zeta\oplus \eta) \in \text{Ker}(L(\jsig,J;v))$, and
$g:\Sig \ra \Sig$ is an element of the automorphism group
$\Group=\text{Aut}(\Sig,\jsig)$, then $\theta_g=g^*\theta$ is also
an element in the kernel of $L$.

Define $\Orbit_\theta$ to be the subspace of  $\Ker(L)$ generated
by $\{\theta_g\}_{g\in \Group}$, and denote by $\m_\theta$ the
ideal of the group ring $\R_\Group$ consisting of all elements

$$
\alpha=\sum_{g \in \mathfrak{G}} a_g.g^{-1} \in \R_\Group,
\ \ \ \ \ \ a_g \in \R,$$
such that
\begin{displaymath}
\sum_{g \in \mathfrak{G}} a_g \theta_g=0.
\end{displaymath}
It is easy to check that $\m_\theta$ is a left ideal of the group
ring $\R_{\mathfrak{G}}$.\\

In this part we will have a discussion of the maximal ideals of $\R_\Group$ and a decomposition
of $\Ker(L)$ into the orbits $\Orbit_\thet$, with  associated ideal being maximal, for the case where
$\Group$ is $\Z_p=\frac{\Z}{p\Z}$, and $p$ is prime.\\

Suppose that $p\neq 2$ is prime. The group $\Z_p$ may be
identified with the multiplicative group of the elements

$$\{1,\lambda,\lambda^2,...,\lambda^{p-1}\}$$ where
$\lambda=e^{\frac{2\pi i}{p}}$. Then the formal sums

$$a_k=\sum_{i=0}^{p-1}\text{Re}(\lambda^{ki}).\lambda^i,
\ \ \ \ \ k=0,1,...,\frac{p-1}{2}$$
and
$$b_k=\sum_{i=0}^{p-1}\text{Im}(\lambda^{ki}).\lambda^i,
\ \ \ \ \ k=1,2,...,\frac{p-1}{2}$$ generate $\R_\Group$ as a
vector space over $\R$. Moreover if $I_k=\langle
a_k,b_k\rangle_\R$ is the subspace generated by $a_k,b_k$, then in
the decomposition $$\R_\Group=I_0\oplus I_1 \oplus ...\oplus
I_{\frac{p-1}{2}},$$ the space of elements which are zero in the
$k$-th component forms an ideal $\m^k$ of $\R_\Group$, which is in
fact
maximal.\\

For any section $\theta \in \Ker(L)$, the sections
$$\theta^k=\sum_{i=0}^{p-1}\text{Re}(\lambda^{ki}).\theta_{\lambda^{i}}$$
and
$$\thet^k=\sum_{i=0}^{p-1}\text{Im}(\lambda^{ki}).\theta_{\lambda^{i}}$$
have the associated ideal $\m^k$.

It is interesting to note that if $\theta^k$ and $\thet^k$ are
zero for $k=1,...,\frac{p-1}{2}$, then
$\m_\theta=\m^0$.\\

Before going further, define a \emph{local product} on $\Ker(L)$,
associated to a point $p\in C$, as follows: Let $\{p_g\}_{g\in
\Group}$ denote the points in the pre-image $\pi^{-1}(p)$ of $p$,
such that $h(p_g)=p_{hg}$. Then define

$$\langle \theta,\thet \rangle_p=
\sum_{g\in \Group}\langle \theta(p_g),\thet (p_g) \rangle,$$
using some metric $\langle .,.\rangle$ on the bundle $E$.\\

Now start from a section $\theta$ in $\Ker(L)$. By changing
$\theta$ with one of $\theta^k$ or $\thet^k$  if necessary, we may
assume that the associated ideal $\m_1=\m_\theta$ is one of $\m^i,
i=0,...,\frac{p-1}{2}$. Call this section $\thet_1$. Choose a
point $p_1$ such that $\thet_1$ is not identically zero at
$(p_1)_g$'s. Then look at the orthogonal complement of
$\Orbit_{\thet_1}$ in $\Ker(L)$, with respect to $\langle
.,.\rangle_{p_1}$. Choose another section $\thet_2$ with the
corresponding ideal $\m_2$ among $\m^i$'s in this orthogonal
complement, together with a point $p_2$ in $C$, different from
$p_1$. Furthermore, choose $p_2$ such that $\thet_2$ is nonzero
above $\pi^{-1}(p_2)$. Then look at the intersection of orthogonal
complements of $\Orbit_{\thet_i}$ with respect to $\langle
.,.\rangle_{p_i}$, $i=1,2$, etc.. This process will give us a
decomposition

\begin{equation}
\begin{split}
&\Ker(L)=\Orbit_{\thet_1}\oplus \Orbit_{\thet_2}
\oplus ... \oplus \Orbit_{\thet_\ell}\\
&\m_{\thet_1}=\m_1,\m_{\thet_2}=\m_2,...,\m_{\thet_\ell}=\m_\ell,
\end{split}
\end{equation}
together with $\ell$ points $p_1,p_2,...,p_\ell$. Note that the
ideals $\m_i$ appearing in this decomposition (the whole
collection) is independent of
the way we decompose.\\

Fix the values $(\jsig,J;v)$ as above such that $L=L(\jsig,J;v)$ has an element
in its kernel. Note that there is an action
of $\Group$ on the spaces
$$\Gamma(\Sig,E)\oplus H^1(\Sig,T\Sig)$$
and
$$\Gamma(\Sig,\Omega^{0,1}_{(\Sig,\jsig)}\otimes_J E).$$

One may also consider the space of those sections which have the
associated ideal equal to a given ideal $\m$ (denoted by
$\Gamma_\m(\bullet)$). It is easy to see that the operator $L$
restricts to an operator $L_\m$,

$$
L_\m=L_\m(\jsig,J;v):\Gamma_\m(\Sig,E)\oplus H_\m^1(\Sig,T\Sig)
\lra \Gamma_\m(\Sig,\Omega_{(\Sig,\jsig)}^{0,1} \otimes_J E)
$$
on the space of such "$\m$-sections". Similarly the adjoint
operator $L^*$ restricts to give an adjoint operator
$L_\m^*$ which goes in the other direction.\\

The index computation of the Appendix shows that for any ideal
$\m$, the operator $L_\m$ is Fredholm of index zero. This fact may
be used to show that, parallel to the above process of decomposing
$\Ker(L)$ into orbits $\Orbit_{\thet_i}$, one may also find
elements $\muu_1,...,\muu_\ell$ in $\Cok(L)$ with the property
that

\begin{equation}
\begin{split}
&\Cok(L)=\Orbit_{\muu_1}\oplus \Orbit_{\muu_2}
\oplus ... \oplus \Orbit_{\muu_\ell}\\
&\m_{\muu_1}=\m_1,\m_{\muu_2}=\m_2,...,\m_{\muu_\ell}=\m_\ell
\end{split}
\end{equation}
with a similar orthogonality assumption (using the same set of
points $p_1,...,p_\ell$).\\

The goal is to show that the space $\DivG_{\{\m_1,...,\m_\ell\}}$
of the tuples $(\jsig,J;v)$ such that the kernel
$\Ker(L(\jsig,J;v))$ has a decomposition as above, is locally a
submanifold of $\mathcal{G}^\Group$. \\

\begin{thm}\label{thm:model}
For odd prime numbers $p$, the space
$\DivG_{\{\m_1,...,\m_\ell\}}$ is an analytic submanifold of
$\mathcal{G}^\Group$, of codimension $2\ell-q$, if $q$ of $\m_i$'s
are equal to $\m^0$.
\end{thm}

\begin{proof}
We give the proof for the case where $\ell=1$ and $\m_1\neq \m^0$.
Then we will make a remark on the other cases. Suppose that
$\Ker(L)=\Orbit_{\thet}$ is generated by $\theta_1,\theta_2$. This
implies that $ \Cok(L)=\Orbit_\muu$ is generated by some
$\mu_1,\mu_2$. We will assume that $\theta_i$'s are orthogonal to
each other and that the same is true for $\mu_i$'s. Furthermore,
assume that $\m_\thet=\m_\muu=\m$.\\

Suppose that $(\jsig',J';v')$ is another element in $\DivGm$ which
is very close to $(\jsig,J;v)$.  Assume that
$\theta_k'=\theta_k+\theta_k^0$ are the corresponding sections of
$\Gamma(\Sig,E)\oplus H^1(\Sig,T\Sig)$. We may assume that
$\theta_k^0$'s are orthogonal to $\Ker(L)$, hence to all
$\theta_i$'s, $i,k\in \{1,2\}$. One may write
$\jsig'=\jsig+\delta+F_1(\delta)$, $J'=J+Y+F_2(Y)$ and
$v'=v+Z+F_3(\delta,Y,Z)$. Then we should have the extra condition that

\begin{equation}
\label{eq:tangentcondition}
\begin{split}
&\delta \jsig +\jsig \delta =0\\
&YJ+JY=0\\
&JZ+Z\jsig+Yv+v\delta=0.
\end{split}
\end{equation}
The functions $F_1,F_2,F_3$ will be analytic functions of their
variables. The vanishing of $L(\jsig',J';v')$ at
$$\theta_k'=(\zeta_k'=\zeta_k+\zeta_k^0)\oplus
(\eta_k'=\eta_k+\eta_k^0)$$ can be written as

\begin{displaymath}
\begin{split}
0=L'(\theta_k')(\omega)=L(\theta_k^0)&(\omega)+
Y\nabla_{\jsig \omega}\zeta_k+
Y\circ d\pi \circ \eta_k\\
&+\frac{1}{2}\{(\nabla_{\zeta_k} Y)(d\pi \circ \jsig) +
(\nabla_{\zeta_k} J)(d\pi \circ \delta)
-\nabla_{\theta_k} Z\}\\
&\ \ \ \ \ \ \ \ \ +\text{terms of higher order},\ \ \ \ k=1,2,
\end{split}
\end{displaymath}
where $L'=L(\jsig',J';v')$ and $L$ is as before.\\

In particular, we have to solve two equations of the form
$L(\theta_k^0)=\bullet$. This is not something that one can  have
any hope to do in general, since $L$ is not surjective. However,
let us follow Taubes (\cite{Taubes}) and denote by $\Pi$ the
projection

$$\Gamma(\Sig,\Omega^{0,1}_{(\Sig,\jsig)}\otimes_J E)
\longrightarrow \Cok(L),$$ and denote the projection over the
image of $L$ by $\Pi^c$. Then $\text{Id}=\Pi+\Pi^c$ and we obtain
the following two equations:

\begin{equation}
\begin{split}
&L(\theta_k^0)+\Pi^c\{Y\nabla_{\jsig \omega}\zeta_k +
Y\circ d\pi \circ \eta_k\\
&\ \ \ \ \ \ \ \ \ \ \ \
+\frac{1}{2}[(\nabla_{\zeta_k} Y)(d\pi \circ \jsig) +
(\nabla_{\zeta_k} J)(d\pi \circ \delta)-\nabla_{\theta_k} Z]\\
&\ \ \ \ \ \ \ \ \ \ \ \ \ \ \ +\text{higher order}\}=0\\
&\Pi \{ Y\nabla_{\jsig \omega}\zeta_k+Y\circ d\pi \circ \eta_k\\
&\ \ \ \ \ \ \ \ \ \ \ \
+\frac{1}{2}[(\nabla_{\zeta_k} Y)(d\pi \circ \jsig) +
(\nabla_{\zeta_k} J)(d\pi \circ \delta)-\nabla_{\theta_k} Z]\\
&\ \ \ \ \ \ \ \ \ \ \ \ \ \ \ +\text{higher order}\}=0\\
\end{split}
\end{equation}

For given values for $(\jsig',J';v')$, or rather for given values
for $(\delta,Y;Z)$, the first equation may be solved uniquely to
give a value for $ \theta_i^0$ as an analytic function of
$(\delta,Y;Z)$. The second set of equations may  then be thought
of as an analytic function

$$F: \mathcal{B}(\jsig,J;v)\subset \mathcal{G}^\Group
\lra \text{Hom}(\Ker(L),\Cok(L))$$
where $ \mathcal{B}(\jsig,J;v)$ denotes a small neighborhood of
$(\jsig,J;v) $ in $\mathcal{G}^\Group$.\\

The map $F$ will be an analytic function and its derivative at
zero will be given explicitly as follows: It will take a tangent
vector $(\delta,Y;Z)$ satisfying the equation
(\ref{eq:tangentcondition}) and will give the matrix with $(i,j)$
entry equal to:

\begin{displaymath}
\begin{split}
[dF]_{ij}=\int_\Sig \langle &\gamma_i,\mu_j  \rangle,\\
\gamma_i= Y\nabla_{\jsig \omega}\zeta_i&+Y\circ d\pi \circ \eta_i \\
&+\frac{1}{2}[(\nabla_{\zeta_i} Y)(d\pi \circ \jsig) +
(\nabla_{\zeta_i} J)(d\pi \circ \delta)-\nabla_{\theta_i} Z].
\end{split}
\end{displaymath}

We will choose $\delta$ and $Y$ to be zero and $Z$ will be any
arbitrary section satisfying $JZ+Z\jsig=0$. Clearly, if such
tuples give us the required surjectivity for $dF$, then
$F^{-1}(0)$ will be a subvariety. More coherently, we may think of
$dF(Z)=dF(0,0;Z)$ as a constant multiple of the pairing
$$\rho_Z:\Ker(L)\times \Cok(L) \lra \R$$
defined by
$$\rho_Z(\theta,\mu)=\int_\Sig \langle \nabla_\theta Z,\mu \rangle.$$
It is important to note that
$\rho_Z(g^*\theta,g^*\mu)=\rho_Z(\theta,\mu)$ for any element $g$
of $\Group$. This implies that  the only independent relations
given by $\rho_Z=0$ are
\begin{displaymath}
\rho_Z(\theta_1,\mu_i)=0,\ \ \ \ \ \ \ \ \ i=1,2.
\end{displaymath}

The equation  $F=0$ is implied from $Q(F)=0$ where $Q$ is the
projection on the first row of $\text{Hom}(\Ker(L),\Cok(L))$. The
differential of $Q(F)$ may be described on $(0,0;Z)$ as the map
taking $Z$ to
$$(\rho_Z(\theta_1,\mu_1),\rho_Z(\theta_1,\mu_2)).$$
We will show that this map is surjective. As a result, $(Q\circ
F)^{-1}(0)=F^{-1}(0)$ is locally an analytic submanifold of
codimension $2$ near the point $(\jsig,J;v)$ of
$\mathcal{G}^\Group$. This is equivalent to showing that the map
taking $Z$ to $(\rho_Z(\thet,\muu_g))_{g\in \Group}$ has rank
$2$.\\

Denote $\nabla_\thet Z$ by $G(\thet)$, where $G$ is a function
depending on $Z$. For a generic choice of $p\in C$, let
$\{p_g\}_{g\in \Group}$ be the points in $\tau_\Sig^{-1}(p)$. Fix
an identification of the fibers $F_p\cong E_{p_g}$ with $\C^3$
such that $J$ becomes the standard complex structure, and an
identification of $T_{p_g}\Sig\cong T_pC$ with $\C$, such that
$\jsig$ is identified with $i$. Denote the $i$-th component of the
image of $\muu(p_g)(1)\in \C^3\cong \R^6$ by
$\muu^i(g)$, $g\in \Group,i=1,...,6$.\\

Look at the values for $Z$ such that they are supported over the
point $p$ (i.e. invariantly supported on $p_g$'s). The
corresponding section $G(\thet)$ will be supported on $p_g$'s as
well. We may assume that

$$\int_\Sig \langle G(\thet),\muu_g\rangle=
\sum_{h\in \Group}f_\thet(p_h)\muu^i(gh),$$ where
$f_\thet:\{p_g\}\ra \R$ is a function that depends on $Z$. Our
freedom in choosing $Z$ guarantees that any such function may be
obtained, except if there is a relation
$$\sum_{g\in \Group} a_g\thet(p_g)=0$$
for the section $\thet$, when we have the similar relation
$$\sum_{g\in \Group} a_g f_\thet(p_g)=0.$$

If there is a relation between $\rho_Z(\thet,\muu_g)$ of the form
$$\sum_{g\in \Group} a_g \rho_Z(\thet,\muu_g)=0,
\ \ \ \ \ \ \ \ \ \ \forall Z$$
then in particular we get

\begin{displaymath}
0=\sum_{g\in \Group}\sum_{h\in \Group}a_g f_\thet(p_h)\muu^i(gh)=
\sum_{h\in \Group} b_h f_\thet(p_h).
\end{displaymath}
Thus, it is implied that $\sum_{h\in \Group}b_h \thet(p_h)=0$.\\

\begin{lem}
If $p$ (and consequently the collection $\{p_g\}_{g\in \Group}$)
come from a generic choice, and $\sum_{h\in \Group}b_h
\thet(p_h)=0$ for $\beta=\sum_{h\in \Group}b_h.h^{-1}$ as above,
then $\beta \in \m$.
\end{lem}
\begin{proof}~(of the lemma).
Suppose that $\beta$ is not in the ideal $\m$, which has a rank
(as a vector space over $\R$) equal to $p-2=|\Group|-2$. Then
there are $(p-1)$ independent relations between the vectors
$\{\thet(p_h)\}_{h\in \Group}$. This implies that $\thet(p_h)$ are
mutually linearly dependent. Without loss of generality we may
assume that $\thet(p_h)=a(p_h).\thet(p_{e})$ where $e\in \Group$
is the identity element. Note that if the above claim is not true
for generic $p$, it will not be true for any $p$. As a result, $a$
defines a function

$$a:\Sig \lra \R, \ \ \ \ \ \ \ \ \thet_h=a_h. \thet=h^*(a). \thet.$$
On the other hand, the above relation implies that

$$0=L(\thet_h)(\omega)=L(a_h.\thet)(\omega)=
a_h L(\thet)+(\omega.a_h)\thet+[(\jsig \omega).a_h](J\thet).$$ As
a result $(\omega.a_h)\thet+[(\jsig \omega).a_h](J\thet)=0$ which
can not be the case for $\thet \neq 0$ and real valued function
$a_h$, unless $a_h$ is constant. If $a_h$ is constant, the
assumptions $p\neq 2$ and $a_h\in \R$ imply that $a_h=1$. Thus
$\thet$ is invariant under the action of $\Group$, contradicting
our assumption. This completes the proof of the lemma.
\end{proof}

We conclude
that $\sum_{h\in \Group}b_h.h^{-1} \in \m_\thet=\m$. But on the other hand

\begin{displaymath}
\sum_{h\in \Group} b_h.h^{-1}=
(\sum_{h\in \Group} \muu^i(h).h^{-1})(\sum_{h\in \Group} a_h.h)
\end{displaymath}
Note that this is true for all components of $\muu$ corresponding
to $i=1,..,6$.   By considering the fact that $\m_\muu=\m$ as
well, the only possible way for

$$(\sum_{h\in \Group} \muu^i(h).h^{-1})(\sum_{h\in \Group} a_h.h)$$
to be in $\m$ for a generic choice of $p$, is when $\sum_{h\in
\Group} a_h.h \in \m$. This proves that the only relations between
$\rho_Z(\thet,\muu_g)$'s are those corresponding to the elements
of $\m$, and completes the proof
of the theorem for the case $\ell=1$.\\

When $\ell>1$, but the ideals are different from $\m^0$, the proof
is just slightly more complicated. In the definition of
orthogonality, one should choose the points $p_i$ (and
correspondingly $(p_i)_g$'s) to be generic for $\thet_i,\muu_i$,
in the sense of the above lemma. Then the independence of the
equations
\begin{displaymath}
\rho_Z(\theta^i_1,\mu^i_1)=0,\ \ \ \
\rho_Z(\theta^i_1,\mu^i_2)=0, \ \ \ \ i=1,2,...,\ell.
\end{displaymath}
follows from the orthogonality, and a discussion similar to the above one.\\

For the ideal $\m^0$ of $\R_\Group$, the corresponding sections
$\thet,\muu$ are the generators of their orbits as a $\R$-vector
space. Here the claim is that $\rho_Z(\thet,\muu)$ is not
identically zero. Similar to the above arguments, if
$\rho_Z(\thet,\muu)=0$ for all $Z$, then at a point $p\in C$, and
the corresponding points $\{p_g\}_{g\in \Group}$, we will have

$$\sum_{g\in \Group}\thet(p_g).\muu^i(p_g)=0$$
where $\muu^i$ denotes the $i$-th component of $\muu$. If
$\sum_{g\in \Group} \muu^i(p_g).g^{-1}$ is not in $\m=\m^0$, then
there are $p$ independent linear relations between
$\{\thet(p_g)\}_{g\in \Group}$, which implies that they are all
zero. Since this can not be true for generic $p$, $$\sum_{g\in
\Group} \muu^i(p_g).g^{-1}\in \m.$$

This is also a contradiction, since it implies that $\sum_{g\in
\Group} \muu^i(p_g).g^{-1}$ is annihilated by all elements of
$\R_\Group$ (it is already annihilated by elements of $\m$, and
being in $\m$ gives one more relation independent of the previous
ones). The passage from one ideal to the $\ell>1$ case is similar
to the above discussion.
\end{proof}

%%%%%%%%%%%%%%%%%%%%%%%%%%%%%%%%%%%
%%%%%%%%%%%%%%%%%%%%%%%%%%%%%%%%%%%
%%%%%%%%%%%%%%%%%%%%%%%%%%%%%%%%%%%

\section{Finiteness}

The goal of this section is to show that for a generic choice of
an almost complex structure $J$ on $X$ and a perturbation term
$v\in \Homi_J(T\PN,TX)$, the somewhere injective solutions to the
equations

\begin{equation}
\begin{split}
&f:(\Sig, \jsig) \lra X,\\
&\pi_\Sig:(\Sig,\jsig)\ra \Ci_g \subset \PN, \\
&f_*[\Sig]=\beta \in H_2(X,\Z),\\
&\dbar_{\jsig,J}f=(\pi_\Sig\times f)^*v,
\end{split}
\end{equation}
are isolated and finite. Our discussion on transversality will
assign a sign to each such solution coming from the transverse
intersection of the moduli spaces and the spectral flow to a
complex $\dbar$ operator. For a more careful treatment of signs,
we refer the reader to \cite{R-T1,R-T2}.\\

Write the homology class $\beta \in P_2(X,\Z)$ as $p\alpha$, where
$\alpha$ is a primitive homology class in $H_2(X,\Z)$ (i.e.
$\alpha$ is not a multiple of some other class). Note that when
$\beta$ itself is primitive, then the standard arguments in
Gromov-Witten theory rules out the possibility of  convergence of
a sequence of embedded solutions to any solution with singular
domain, and we can not have a multiply covered solution. So the
set of embedded solutions for a generic choice of parameters
$(J,v)$ consists of isolated points and is compact, thus finite.
The independence of the signed count of these solutions from the
parameters is also standard.
So we will assume that $p\neq 2$ is a prime number.\\

Suppose that $(J,v,(f:(\Sig,\jsig)\xra{\tau_\Sig}C\xra{h}X))$ is
an element in $\ModGX$. We will get a map from a neighborhood of

$$(J,v,(f:(\Sig,\jsig)\xra{\tau_\Sig}C\xra{h}X))$$ in $\ModGX$ to
the local model $\mathcal{G}^\Group$ associated with the bundle
$E=f^*TX$. The map takes this point to $(\jsig,\bar J,\bar v)$,
where $\bar J$ is the induced almost complex structure on $f^*TX$
and $\bar v$ is the perturbation term induced by $v$ on a an
identification of $h^*TX$ with a tubular neighborhood of $h(C)$.
Denote this map by

$$q:\mathcal{B}(J,v,(f,\jsig)) \subset\ModGX \lra \mathcal{B}
(\jsig,\bar J,\bar v) \subset  \mathcal{G}^\Group. $$ Here
$\Group$ is the automorphism group  of $(\Sig,\jsig)$. Construct a
map $F$ from $\mathcal{B}(\jsig,\bar J,\bar v)$ to
$\text{Hom}(\Ker(L),\Cok(L))$ as in the proof of
theorem~\ref{thm:model}, where $L=L(\jsig,\bar J,\bar v)$. Also
construct the projection $Q$ as in the proof of
theorem~\ref{thm:model}. Since in the proof of surjectivity of
$d_{(\jsig,\bar J,\bar v)}(Q\circ F)$, we only used the
perturbation of $\bar v$, it can be easily
concluded that the composition $Q\circ F \circ q$ has a surjective
derivative at $(J,v,(f,\jsig))$.\\

This observation implies that we obtain submanifolds

$$\overline{\mathcal{D}}_{\{\m_1,...,\m_\ell\}}^\Group \subset \ModGX$$
consisting of the points $(J,v,(f,\jsig))$ such that the
linearization operator

\begin{displaymath}
\begin{split}
L=L(J,v,(f,\jsig)):H^1_{\jsig}(\Sig,T\Sig)\oplus \Gamma&(\Sig,f^*TX)\\
&\lra
\Gamma(\Sig,\Omega_{(\Sig,\jsig)}^{0,1}\otimes_Jf^*TX)
\end{split}
\end{displaymath}
defined by

\begin{equation}
\begin{split}
L(\zeta,\eta)(\omega)=\nabla_\omega \zeta&+
J\nabla_{\jsig \omega}\zeta+(J\circ df \circ \eta) (\omega)\\
&+\frac{1}{2}\{(\nabla_\zeta J)(df \circ \jsig)-
(\nabla_{(\zeta \oplus \eta)}v)\},
\end{split}
\end{equation}
has a kernel $\Ker(L)=\Orbit_{\thet_1}\oplus ...\Orbit_{\thet\ell}$.
Here we assume that $\m_{\thet_i}=\m_i$.\\

The submanifold
$\overline{\mathcal{D}}_{\{\m_1,...,\m_\ell\}}^\Group$ will have a
codimension equal to
$\text{dim}(\Ker(L))$, using theorem~\ref{thm:model}.\\

After setting up the above notation, we are now ready to prove the
following theorem (compare with theorem~\ref{thm:definition}),
which shows that the claimed counts of the embedded solutions of
the perturbed Cauchy-Riemann equation, are in fact meaningful.

\begin{thm}
\label{thm:generic}
For a Bair subset $\ParMod_{reg}\subset
\ParMod$, the following is true: If $(J;v)\in \ParMod_{reg}$, then
the space $\ModXpar$ of the somewhere injective solutions
$$f:(\Sig,\jsig)\lra X$$ to the perturbed Cauchy-Riemann equation
$$\dbar_{\jsig,J}f=(\tau_\Sig \times f)^*v, \ \ \ \ \ \ \
f_*[\Sig]=\beta, $$ is finite. Moreover, at any such solution,
$L=L(J,v,(f,\jsig))$ has a trivial kernel. Finally, to any such
solution $(J,v,(f,\jsig))$ is assigned a sign
$\epsilon(f,\jsig)=\epsilon(J,v,(f,\jsig))$ coming from the
spectral flow from $L$ to a complex $\dbar$ operator.
\end{thm}
\begin{proof}
Consider all of the possible actions of the group
$\Z_p=\frac{\Z}{p\Z}$ on the surface $\Sigma$. Denote such an
action by $\Group$ and fix the quotient map $\pi:\Sig \ra
\frac{\Sig}{\Group}$. There are finitely many such group actions.
For any such action $\Group$, consider the manifold $\ModGX$ and
the projection map
$$q_{\Group}:\ModGX \lra \ParMod.$$
It is easy to check that $q_{\Group}$ is in fact a Fredholm operator of index zero.\\

For any set $\{\m_1,...,\m_\ell\}$ of ideals, define
$d(\{\m_1,...,\m_\ell\})$ to be the dimension of the corresponding
kernel. Then if we restrict the map $q_\Group$ to the submanifold
$\mathcal{D}=\overline{\mathcal{D}}^\Group_{\{\m_1,...,\m_\ell\}}$,
the index of this restriction will be equal to
$-d(\{\m_1,...,\m_\ell\})$. The set of regular values for all the
projection maps $q_\Group$ and $q_\Group|_{\mathcal{D}}$, for
different choices of the group action $\Group$ and the ideals
$\{\m_1,...,\m_\ell\}$, will still be a Bair subset
$\ParMod_{reg}\subset \ParMod$ (note that in particular we
consider the case where $\Group$ is the trivial group). If $(J,v)$
is a regular value of $q_\Group|_{\mathcal{D}}$, then
$q_\Group^{-1}(J,v)\cap {\mathcal{D}}=\emptyset$, since the index
of
$q_\Group|_{\mathcal{D}}$ is negative.\\

Suppose that $(J;v)\in \ParMod_{reg}$. Then the points in
$\ModXpar$ will be isolated, and at any such point, the kernel of
the linearization map is trivial. We will be done if we can show
that $\ModXpar$ has only finitely many points, since the signs
$\epsilon(J,v,(f,\jsig))$ may be assigned as in
\cite{R-T1,R-T2}.\\

Suppose that this is not the case and $\{f_n:(\Sig,{\jsig}_n)\ra
X\}$ is a sequence of somewhere injective solutions to the
perturbed Cauchy-Riemann equation above. This sequence will have a
convergent subsequence. The limit will be a map

$$f=f_\infty:(\overline{\Sig},\jsig) \lra X$$ which solves the
same equation. We have already argued that $\overline{\Sig}$ can
not be singular and we may identify it with $\Sig$.\\

If $f$ is somewhere injective then $(f,\jsig)\in \ModXpar$. By
lemma~\ref{lem:uniqueness}, $L(J,v,(f,\jsig))$
has to have a nontrivial kernel, which contradicts our assumption on $(J;v)$.\\

So, $f$ will factor as

$$f:\Sig \xra{\tau} \frac{\Sig}{\Group} \xra{g}, X$$
where $\Group$ denotes a group action on $(\Sig,\jsig)$. But
$\beta=p \alpha$, and $\alpha$ is primitive. This implies that the
underlying group of  $\Group$ is $\Z_p=\frac{\Z}{p\Z}$. So
$(J,v,(f,\jsig))\in \ModGX$. Again the convergence of a sequence
of solutions to $(f,\jsig)$ implies that the kernel of the
linearized operator is nontrivial. Thus

$$(J,v,(f,\jsig))\in \overline{\mathcal{D}}^\Group_{\{\m_1,...,\m_\ell\}}$$
for some choice of ${\{\m_1,...,\m_\ell\}}$. This is a
contradiction, completing the proof of the theorem.
\end{proof}

Using the information given by the above theorem, we define:
\begin{defn}
For a prime homology class $\beta \in P_2(X,\Z)$, choose a point
$(J,v)\in \ParMod_{reg}$. Define
$$\Inv_g(\beta)=\sum_{(f,\jsig)\in \ModXpar}\epsilon(f,\jsig)$$
\end{defn}

%%%%%%%%%%%%%%%%%%%%%%%%%%%%%%%%%%%
%%%%%%%%%%%%%%%%%%%%%%%%%%%%%%%%%%%
%%%%%%%%%%%%%%%%%%%%%%%%%%%%%%%%%%%
\section{Invariance}

This section will be devoted to the proof of the invariance of
$$\Inv_g:P_2(X,\Z) \lra \Z$$ from the choice of the regular values
 $(J,v)\in \ParMod_{reg}$. We begin this section
by setting up a way of thinking of the moduli space of solutions
corresponding to paths between two regular values
$\gamma_i=(J_i,v_i)\in \ParMod_{reg}$ for $i=0,1$. Namely, denote
by $\PathMod$ the moduli space of the paths

$$\gamma:[0,1]\lra \ParMod,\ \ \ \ \text{with} \
\gamma(i)=\gamma_i,\  i=0,1.$$
There is a section

$$\dbar:\mathcal{Z}=[0,1] \times \PathMod
\times \Mod_g  \times \MapMod \lra \Bundle$$ defined by
$\dbar(t,\gamma,\jsig,f)=\dbar(\gamma(t),\jsig,f)$. Here we will
think of $\Bundle$ as the bundle pulled back to $\mathcal{Z}$ from
$\mathcal{Y}$, via the map $\mathcal{Z}\ra \mathcal{Y}$ defined by

$$(t,\gamma,\jsig,f) \ra (\gamma(t),\jsig,f).$$

Again, we may also define the equivariant versions of these,
giving the maps

$$\dbar_\Group
:\mathcal{Z}^\Group=[0,1] \times \PathMod  \times \Mod^\Group
\times \MapMod^\Group \lra \Bundle^\Group.$$

One may argue, using an argument similar to those in the section
on transversality, that the sections $\dbar$ and $\dbar_\Group$
are transverse to the zero section. As a result, we obtain a
smooth submanifold of $\mathcal{Z}$, consisting of the zeros of
$\dbar$, which will be denoted by $\NodX$. Similarly we may define
$\NodGX$
as a submanifold of $\mathcal{Z}^\Group$.\\

Note that theorem~\ref{thm:invariance} is a direct corollary of
the following theorem:

\begin{thm}
\label{thm:invariance2}
There is a Bair subset
$$\PathMod_{reg}\subset \PathMod,$$ such that for $\gamma \in
\PathMod_{reg}$, the moduli space $\NodXpar$ of the somewhere
injective solutions associated with the path
$\gamma(t)=(J_t,v_t)$, forms a compact $1$-manifold, giving a
cobordism between $\mathcal{M}_g(X,\beta;(J_0,v_0))$ and
$\mathcal{M}_g(X,\beta;(J_1,v_1))$. In particular $\Inv_g(\beta)$,
as computed using $(J_0,v_0)$ is equal to $\Inv_g(\beta)$, as
computed using $(J_1,v_1)$.
\end{thm}

\begin{proof}
There are projection maps from $\NodX$ to $\PathMod$, and from
$\NodGX$ to $\PathMod$, which are Fredholm maps of index $1$. As a
result there is a Bair set of regular values of these projection
maps which we will denote by $\PathMod_{reg}$. The pre-image of
our submanifolds $\Divl$ will be submanifolds of $\NodGX$, and the
restriction of the projection map will typically have negative
index on these submanifolds. The exception is
$\overline{\mathcal{D}}_{\m^0}^{\Group}$, where the index is zero.
We may choose the set of regular values so that
they are regular values of these restriction maps as well.\\

With this choice of $\PathMod_{reg}$, if $\gamma \in
\PathMod_{reg}$ then $\NodGXpar$ will be a smooth $1$-dimensional
manifold and for any $(t,\gamma,\jsig,f) \in \NodGXpar$ the
linearization $L(\gamma(t),(f,\jsig))$ will have a trivial kernel,
except for an isolated set of such points where the kernel is
$1$-dimensional.\\

We are interested in a study of a neighborhood of the points
$(t,\gamma,\jsig,f)$ giving $(J,v;(f,\jsig))$ (with
$(J,v)=\gamma(t)$) where the kernel of $L=L(J,v,(\jsig,f))$ is
one dimensional.\\

Suppose that this kernel is generated by $\theta=(\zeta,\eta)$ and
that the corresponding cokernel is denoted by $\mu$. At the time
$t+\epsilon$, a nearby point corresponding to the parameter
$\gamma(t+\e)=(J_\e,v_\e)$ may be described as follows: If $(Y,Z)$
denotes the derivative $\gamma'(t)$,
then $J_\e=J+\e Y +G_1(\e), v_\e=v+\e Z +G_2(\e)$.
Here $G_1(\e)$ and $G_2(\e)$ are in $o(|\e|)$.\\
Any solution to the perturbed Cauchy-Riemann equation
corresponding to $(J_\e,v_e)$ which is close to $(f,\jsig)$ will
be of the form $(f_\e,j_\e)=\text{exp}_{(f,\jsig)}(\thet)$.
$\thet$ may be written as $\thet=s \theta + \theta_0$, where
$\theta_0=(\zeta_0,\eta_0)$ is orthogonal to the kernel of $L$.
The equation

$$\dbar(J_\e,v_\e,(f_\e,j_\e))=0$$
may be reformulated as

\begin{displaymath}
L(s\theta + \theta_0)+\e(Y\circ df \circ \jsig -
(\tau_\Sig \times f)^*Z)+\text{higher order terms}=0
\end{displaymath}

Again let $\Pi,\Pi^c$ be the projection over the kernel of
$L^*$ and image of $L$, respectively. Then we may
rewrite the above equation as the following two equations:

\begin{equation}
\begin{split}
&L(\theta_0)+\Pi^c\{ \e(Y\circ df \circ \jsig -
(\tau_\Sig \times f)^*Z)+\text{higher order terms}\}=0,\\
&\Pi \{\e(Y\circ df \circ \jsig -(\tau_\Sig \times f)^*Z)
+\text{higher order terms} \}=0.
\end{split}
\end{equation}

The first equation gives $\theta_0$ uniquely as an analytic
function of $s,\e$ with no linear terms in $s$. As a result, The
second equation will be an analytic function of $s,\e$ which we
denote by $g(s,\e)$. An argument similar to the one used by Taubes
in \cite{Taubes} shows that the Taylor expansion of $g$ starts as

$$g(s,\e)=r_1 \e + r_2 s^2 +\text{higher order terms}.$$

Note that $r_1$ is in fact computed as

\begin{displaymath}
r_1=\int_\Sig \langle (Y\circ df \circ \jsig -
(\tau_\Sig \times f)^*Z),\mu \rangle
\end{displaymath}
Our assumption on the regularity of $\gamma$ for all projection
maps implies that this pairing is nonzero. An argument similar to
that of \cite{Taubes}, shows that $r_2$ is also proportional to
the differential of the map $Q\circ F$, which is pulled back to
the submanifold $\NodXpar$. The regularity implies again, that
$r_2$ is also nonzero. This gives a description of the
neighborhood of $(t,\gamma,\jsig,f)$
in $\NodXpar$.\\

Note that this discussion may be followed even if
$(t,\gamma,\jsig,f)$ is a weak limit point
of a sequence in $\NodXpar$ that lies in  $\NodGXpar$.\\

If $(t,\gamma,\jsig,f)$ is  in $\NodXpar$, then the projection

$$pr:\NodXpar \lra [0,1]$$ will locally be  like the map going
from $\{(s,\e)|\ r_1 \e+r_2 s^2=0 \}$ to $[0,1]$ which sends
$(s,\e)$ to $t+\e$. As a result $(t,\gamma,\jsig,f)$ is a critical
point of the projection
map $pr$ and nothing interesting happens at this point.\\

Now suppose that $(t,\gamma,\jsig,f)$ is in $\NodGXpar$, with
$\Group=\Z_p$. Since the kernel is nontrivial, $(J,v,(f,\jsig))$
is forced to be in $\Divisor$. Note that the other submanifolds
$\Divl$ are excluded by the regularity. As a result, the sections
$\theta,\mu$ are $\Group$-invariant and the whole argument above
may be done for the $\Group$ invariant sections, and inside
$\NodGXpar$. The result is that a neighborhood of
$(t,\gamma,\jsig,f)$ in $\NodGXpar$ is described by a similar
function $g'(s,\e)=0$. A version of the uniqueness
lemma~\ref{lem:uniqueness} may be used to conclude that this
neighborhood of $(t,\gamma,\jsig,f)$ in $\NodGXpar$ is in fact
identical to the previous neighborhood of it (as a weak limit
point) in $\NodXpar$. This is a contradiction which proves the theorem.
\end{proof}

\begin{remark}
When $p=2$, there is a possibility that at a point
$(t,\gamma,\jsig,f)$ of $\mathcal{N}_g^{\Z_2}(X,\beta;\gamma)$
where the linearization has a kernel with the corresponding ideal
$m^1$, a solution in $\NodXpar$ is blown up for the times $t+\e$
with $\e>0$, while there is no such solution for the times $t-\e$
which is close to $(t,\gamma,\jsig,f)$. This means that to count
the embedded solutions in a class $\beta=2\alpha$, we will always
get contributions from the curves in the class $\alpha$. As a
result a passage from the arguments of this paper to the case
where $p=2$ requires a weighted count of the solutions similar to
the counts in \cite{Taubes}. We postpone such wall-crossing formulas to a
future paper.
\end{remark}

%%%%%%%%%%%%%%%%%%%%%%%%%%%%%%%%%%%
%%%%%%%%%%%%%%%%%%%%%%%%%%%%%%%%%%%
%%%%%%%%%%%%%%%%%%%%%%%%%%%%%%%%%%%
%\input{EmbeddedToGW}
%%%%%%%%%%%%%%%%%%%%%%%%%%%%%%%%%%%
%%%%%%%%%%%%%%%%%%%%%%%%%%%%%%%%%%%
%%%%%%%%%%%%%%%%%%%%%%%%%%%%%%%%%%%
\section{Appendix: A Riemann-Roch formula}

{\large{\bf{Introduction}}}\\

Suppose that $(\Sigma,j=j_\Sigma)$ is a Riemann surface with
automorphism group $\mathfrak{G}$. Denote the quotient curve
$\Sigma/\mathfrak{G}$ by $C$. Fix a divisor $D$ on $\Sigma$ which
is invariant under the action of $\mathfrak{G}$. This means that
for any $\sigma:\Sigma \rightarrow \Sigma$ in
the automorphism group $\mathfrak{G}$ we have $\sigma ^*(D)=D$.\\

For any section $\eta$ of the line bundle $[D]$ over $\Sigma$,
note that the pull back $\eta_\sigma=\sigma^*\eta$ is also a
section of $[D]$. Let $\m_\eta$ be the ideal of the group ring
$\R_{\mathfrak{G}}$ consisting of all elements
\begin{displaymath}
\alpha=\sum_{\sigma \in \mathfrak{G}} a_\sigma.\sigma^{-1},
\ \ \ \ \ \ a_\sigma \in \R,
\end{displaymath}
such that
\begin{displaymath}
\sum_{\sigma \in \mathfrak{G}} a_\sigma \eta_\sigma=0.
\end{displaymath}

It is easy to check that $\m_\eta$ is a left ideal of the
group ring $\R_{\mathfrak{G}}$.\\

Fix an ideal $\mathfrak{m}$ of the group ring
$\mathbb{R}_{\mathfrak{G}}$ and let

\begin{displaymath}
H^i_{\mathfrak{m}}([D]):=\Big\{\eta \in H^i(\Sigma,[D]) \ \Big | \
\mathfrak{m}\subset \m_\eta \Big \},\ \ \ \ i=0,1.
\end{displaymath}

Define
\begin{displaymath}
h^i_{\mathfrak{m}}([D])=\text{dim}_{\mathbb{R}}
(H^i_{\mathfrak{m}}([D])),
\end{displaymath}
and let $\chi_{\mathfrak{m}}([D])=
h^0_{\mathfrak{m}}([D])-h^1_{\mathfrak{m}}([D])$.\\

Our goal in this note is to obtain some results on the behavior of
$\chi_{\mathfrak{m}}([D])$ in terms of the topological properties
of $[D]$ and the covering
\begin{displaymath}
\pi:\Sigma \rightarrow C=\frac{\Sigma}{\mathfrak{G}}.
\end{displaymath}

Let us start with a consideration of the holomorphic tangent
bundle $T_\Sigma$ of $\Sigma$. If $B$ denotes the branching
divisor of the covering map $\pi:\Sigma \rightarrow C$, then we
will have
\begin{displaymath}
T_\Sigma=\pi^*T_C\otimes [-B].
\end{displaymath}
This gives the short exact sequence:
\begin{displaymath}
0 \longrightarrow T_\Sigma \longrightarrow \pi^* T_C
\longrightarrow \mathcal{O}_B \longrightarrow 0.
\end{displaymath}

If we consider the sections corresponding to the left ideal
$\mathfrak{m}$ of $\mathbb{R}_{\mathfrak{G}}$, we will get the
following long exact sequence:

\begin{displaymath}
\begin{split}
0 &\longrightarrow H^0_{\mathfrak{m}}(T_\Sigma)
\longrightarrow H^0_{\mathfrak{m}}(\pi^*T_C) \longrightarrow
H^0_{\mathfrak{m}}(\mathcal{O}_B)\longrightarrow \\
&\longrightarrow H^1_{\mathfrak{m}}(T_\Sigma)
\longrightarrow H^1_{\mathfrak{m}}(\pi^*T_C) \longrightarrow
H^1_{\mathfrak{m}}(\mathcal{O}_B)=0.
\end{split}
\end{displaymath}

As a result
\begin{equation}
\chi_{\mathfrak{m}}(T_\Sigma)=\chi_{\mathfrak{m}}
(\pi^*TC)-h^0_{\mathfrak{m}}(\mathcal{O}_B).
\end{equation}

In the next step we consider line bundles of the form $\pi^*L$
where $L\rightarrow C$ is a line bundle on the quotient curve $C$.
Any such line bundle may be written as $L=[E]$ where $E$ is a
divisor on $C$ such that its support is disjoint from the branched
locus. Note that if $p$ is a point on $C$ which is not in the
branched locus of $\pi:\Sigma \rightarrow C$ then there is a short
exact sequence:

\begin{displaymath}
0\longrightarrow \pi^*L=[\pi^*E]
\longrightarrow [\pi^*(E+p)]\longrightarrow
\mathcal{O}_{\pi^*p}=
\bigotimes_{\sigma \in \mathfrak{G}}\mathcal{O}_{p_\sigma}
\longrightarrow 0,
\end{displaymath}
where $\{p_\sigma\}_{\sigma \in \mathfrak{G}}$ is the
pre-image $\pi^{-1}\{p\}$.\\

Again taking the long exact sequence corresponding to the ideal
$\mathfrak{m}$ gives

\begin{equation}
\chi_{\mathfrak{m}}(\pi^*[E+p])=\chi_{\mathfrak{m}}
(\pi^*[E])+h^0_{\mathfrak{m}}([\pi^*p]).
\end{equation}

Note that both $\mathbb{R}_{\mathfrak{G}}$ and $\mathfrak{m}$ are
vector spaces over $\mathbb{R}$. The dimension of the first one is
$\mathfrak{g}=|\mathfrak{G}|$ and the dimension of the second one,
we denote by $\mathfrak{r}$. It is easy to check that since $p$ is
a generic point,
$h^0_{\mathfrak{m}}([\pi^*p])=2(\mathfrak{g}-\mathfrak{r})$. As a
result of these two observations
\begin{equation}
\chi_{\mathfrak{m}}(\pi^*L)=2(\mathfrak{g}-\mathfrak{r}).\text{deg}(L)
+\chi_{\mathfrak{m}}(\mathcal{O}_\Sigma).
\end{equation}

In the last step, we compare the line bundle $\mathcal{O}_\Sigma$
with the canonical bundle $K_\Sigma$. Note that
$K_\Sigma=\pi^*K_C+[B]$ where $B$ is the branching divisor
introduced earlier. From the short exact sequence
\begin{displaymath}
0 \longrightarrow \pi^*K_C \longrightarrow K_\Sigma
\longrightarrow \mathcal{O}_B \longrightarrow 0,
\end{displaymath}
and the consideration of the global sections, we obtain

\begin{displaymath}
\begin{split}
\chi_{\mathfrak{m}}(K_\Sigma)&
=\chi_{\mathfrak{m}}(\pi^*K_C)+h^0_{\mathfrak{m}}(\mathcal{O}_B)\\
   &=2(\mathfrak{g}-\mathfrak{r}).\text{deg}(K_C)
   +\chi_{\mathfrak{m}}(\mathcal{O}_\Sigma)
   +h^0_{\mathfrak{m}}(\mathcal{O}_B)\\
   &=2(\mathfrak{g}-\mathfrak{r})(2h-2)
   +\chi_{\mathfrak{m}}(\mathcal{O}_\Sigma)
   +h^0_{\mathfrak{m}}(\mathcal{O}_B).
\end{split}
\end{displaymath}

On the other hand, by Serre duality
$\chi_{\mathfrak{m}}(\mathcal{O}_\Sigma)+\chi_{\mathfrak{m}}(K_\Sigma)=0$.
This implies that
\begin{displaymath}
-h^0_{\mathfrak{m}}(\mathcal{O}_B)
=2(\mathfrak{g}-\mathfrak{r})(2h-2)+2\chi_{\mathfrak{m}}(\mathcal{O}_\Sigma)
\end{displaymath}

Combining with the information on $T_\Sigma$ and the fact that
$h^0_{\mathfrak{m}}(T_\Sigma)=0$, this implies that
\begin{equation}
h^1_{\mathfrak{m}}(T_\Sigma)=-3\chi_{\mathfrak{m}}(\mathcal{O}_\Sigma).
\end{equation}

{\large{\bf{General invariant bundles; The index computation}}}\\

Now suppose that $E\rightarrow C$ is a bundle over the quotient
surface, of rank $n$. The bundle $\pi^*E$ will be invariant under
the action of the automorphism group $\mathfrak{G}$ and one may
consider the global sections associated with a left ideal
$\mathfrak{m}$ of the group ring $\mathbb{R}_{\mathfrak{G}}$.
Namely:
\begin{displaymath}
H^i_{\mathfrak{m}}(\Sigma,\pi^*E)
:=\{\eta \in H^i(\Sigma,\pi^*E) \ | \ \mathfrak{m}\subset \m_\eta \}.
\end{displaymath}
The Euler characteristic $\chi_{\mathfrak{m}}(\pi^*E)$
may be defined similarly. \\

For any such bundle, we may formally break it down to the line
bundles:
\begin{displaymath}
E=L_1\oplus L_2\oplus ...\oplus L_n.
\end{displaymath}

From this presentation

\begin{displaymath}
\begin{split}
\chi_{\mathfrak{m}}(\pi^*E)&
=\sum_{i=1}^{n}\chi_{\mathfrak{m}}(\pi^*L_i)\\
&=\sum_{i=1}^{n}[2(\mathfrak{g}-\mathfrak{r}).\text{deg}(L_i)
+\chi_{\mathfrak{m}}(\mathcal{O}_\Sigma)]\\
&=2(\mathfrak{g}-\mathfrak{r}).c_1(E)
+n(\chi_{\mathfrak{m}}(\mathcal{O}_\Sigma))
\end{split}
\end{displaymath}

In the last part of this note we will consider the index
computation associated with the ideal $\mathfrak{m}$
of the group ring $\mathbb{R}_{\Group}$.\\

In our moduli problems, there is a kernel isomorphic to
\begin{displaymath}
H^1_{\mathfrak{m}}(\Sigma,T_\Sigma)\oplus
H^0_{\mathfrak{m}}(\Sigma,\pi^*(TX|_{C})),
\end{displaymath}
where our quotient curve $C$ is embedded in an almost complex
symplectic manifold $(X,\omega,J)$, and
$TX$ is the tangent bundle of $X$.\\

The cokernel will be isomorphic to
\begin{displaymath}
H^0_{\mathfrak{m}}(\Sigma,\pi^*(TX|_{C})).
\end{displaymath}

This implies that the index $I$ of our operator is equal
to

\begin{equation}
\begin{split}
{I}&=h^1_{\mathfrak{m}}(T_\Sigma)+\chi_{\mathfrak{m}}(\pi^*(TX|_{C}))\\
 &=2(c_1(X).[C])(\mathfrak{g}-\mathfrak{r})
 +(n-3)\chi_{\mathfrak{m}}(\mathcal{O}_\Sigma),
\end{split}
\end{equation}
where $[C]$ represents the homology class represented by the curve
$C$ in $H_2(X,\mathbb{Z})$ and
$n$ is the complex dimension of $X$.\\

In particular if $c_1(X)=0$ and $n=3$ then $I=0$:

\begin{cor}
The index of the linearized operator associated with any left
left ideal $\mathfrak{m}$ of the group ring
$\mathbb{R}_{\mathfrak{G}}$ is equal to zero, as far as the target
manifold is a symplectic $3$-fold with vanishing first Chern class
$c_1(X)=0$, the kernel is isomorphic to
 $$H^1_{\mathfrak{m}}(\Sigma,T_\Sigma)\oplus
 H^0_{\mathfrak{m}}(\Sigma,\pi^*(TX|_{C}))$$
 and the cokernel is isomorphic to
$$H^0_{\mathfrak{m}}(\Sigma,\pi^*(TX|_{C})).$$
\end{cor}

%------------------------------------------------------------------
%%%%%%%%%%%%%%%%%%%%%%%%%%%
%\bibliographystyle{amsplain}
%\bibliography{biblio.bib}

\end{document}